\documentclass[12pt,english]{article}
\usepackage[margin=0.75in]{geometry}
\usepackage[utf8]{inputenc}
\usepackage[english]{babel}
\usepackage{amsthm}
\usepackage{graphicx}
\usepackage{url}
\usepackage{hyperref}
\usepackage{amsmath}
\usepackage{listings}
\usepackage{amsfonts}
\usepackage{xcolor}
\usepackage{float}
\newtheorem{theorem}{Theorem}

\def\XXint#1#2#3{{\setbox0=\hbox{$#1{#2#3}{\int}$}
     \vcenter{\hbox{$#2#3$}}\kern-.5\wd0}}

\begin{document}

\title{On the series expansion of k-free Dirichlet series and its analytical continuation}

\author{Artur Kawalec}

\date{}
\maketitle

\begin{abstract}
In this article, we develop a k-free zeta Dirichlet series into a Laurent series with a simple pole, and prove a Stieltjes like formula for the expansion coefficients of the regular part. We also investigate another analytical continuation of these series and develop a formula for $\zeta(\tfrac{1}{k})$ for positive integer $k\geq 2$ in terms of the k-free indicator function.
\end{abstract}

\section{Introduction}
In the previous article [2], we developed several results concerning the square-free Dirichlet series, such as its Laurent series expansion and analytical continuation, and also a formula for $\zeta(\frac{1}{2})$ represented by a certain Dirichlet series with $|\mu(n)|$ coefficients. In this article, we generalize the same results but for the k-free case for positive integer $k\geq 2$. We have the Dirichlet series
\begin{equation}\label{eq:1}
\frac{\zeta(s)}{\zeta(ks)}=\sum_{n=1}^{\infty}\frac{\mu^{(k)}(n)}{n^s}
\end{equation}
which is absolutely convergent for $\Re(s)>1$. Let $\mu(n)$ be the standard M\"obius function. We recall that an integer $n$ is k-free (for $k\geq 2$) if it is not divisible by any $k^{th}$ power equal to $k$ [3][5][6]. We write an indicator function

\begin{equation}\label{eq:1}
\mu^{(k)}(n) = \begin{cases}
  1  &  \text{if n is k-free} \\
  0 &  \text{otherwise}
\end{cases}
\end{equation}
to distinguish between the k-free numbers. The k-free indicator function can be defined in terms of the M\"obius function by the relation
\begin{equation}\label{eq:1}
\mu^{(k)}(n) = \sum_{d^k|n}\mu(d)
\end{equation}
For the square-free integer case $k=2$ one has $\mu^{(2)}(n)=\mu^{2}(n)=|\mu(n)|$. And for $k=3$ case $\mu^{(3)}(n)$ is referred to as cube-free and so forth.

We define the Laurent series expansion

\begin{equation}\label{eq:1}
\frac{\zeta(s)}{\zeta(ks)}=\frac{1}{\zeta(k)(s-1)}+\sum_{n=0}^{\infty}\gamma_n^{M,k} \frac{(s-1)^n}{n!}
\end{equation}
about the pole at $s=1$ and residue $\frac{1}{\zeta(k)}$.  This series is analytic on a disk centered at $s=1$ and radius of convergence $R=1+\frac{2}{k}$ due being limited by the next closest to the pole at the first trivial zero when $ks=-2$ for $s=-\frac{2}{k}$ from the center at $s=1$.

The main result of this article is we derive a general log-limit formula for these expansion coefficients as

\begin{theorem}
\label{thm:fta}
For $n\geq 0$, the coefficients are:
\begin{equation}\label{eq:1}
\gamma^{M,k}_n  =\lim_{x\to\infty}(-1)^n\Bigg\{\sum_{j=1}^{x}\frac{\mu^{(k)}(j)\log^n(j)}{j}-\frac{\log^{n+1}(x)}{\zeta(k)(n+1)}\Bigg\}
\end{equation}

\end{theorem}

\noindent The limit formula (5) was discovered by Wolf in [7] for $k=2$ case, and as we showed in [2]. In this article, we extend the same Wolf formula for the general k-free $k\geq 2$ case.

We now consider the case for $n=0$, where the growth of this Dirichlet series

\begin{equation}\label{eq:1}
\sum_{n\leq x}^{}\frac{\mu^{(k)}(n)}{n}=O(\log x)
\end{equation}
gives a crude log bound; the reason because it's always less than the harmonic series
\begin{equation}\label{eq:1}
\sum_{n\leq x}^{}\frac{\mu^{(k)}(n)}{n}<\sum_{n\leq x}^{}\frac{1}{n}=O(\log x).
\end{equation}
The harmonic series also has a sharper growth estimate
\begin{equation}\label{eq:1}
\sum_{n\leq x}^{}\frac{1}{n}=\gamma+\log x +O(\frac{1}{x})
\end{equation}
where is $\gamma$ is the Euler-Mascheroni constant, and the extraction of this constant is by the limit

\begin{equation}\label{eq:1}
\begin{aligned}
\gamma &=\lim_{x\to\infty} \left(\sum_{n\leq x}^{}\frac{1}{n}-\log x\right) \\
&=0.57721566490153286060\ldots .
\end{aligned}
\end{equation}
In [7] Wolf considered a similar square-free analogue case

\begin{equation}\label{eq:1}
\sum_{n\leq x}^{}\frac{|\mu(n)|}{n}=\gamma^M +\frac{6}{\pi^2}\log x +O(\frac{1}{x})
\end{equation}
where the new constant
\begin{equation}\label{eq:1}
\begin{aligned}
\gamma^M &=\lim_{x\to\infty} \left(\sum_{n\leq x}^{}\frac{|\mu(n)|}{n}-\frac{6}{\pi^2}\log x\right)\\
&=1.04389451571193829740\ldots
\end{aligned}
\end{equation}
is an analogue of the Euler's constant, where it is associated with $k=2$ square-free case. We now extend these results to the general k-free case and write an estimate for the k-free harmonic sum

\begin{equation}\label{eq:1}
\sum_{n\leq x}^{x}\frac{\mu^{(k)}(n)}{n}=\gamma^{M,k}+\frac{1}{\zeta(k)}\log x+O(\frac{1}{x})
\end{equation}
where the generalized k-free analogue of Euler constant $\gamma^{M,k}$ is extracted the same way

\begin{equation}\label{eq:1}
\gamma^{M,k} =\lim_{x\to\infty} \left(\sum_{n=1}^{x}\frac{\mu^{(k)}(n)}{n}-\frac{1}{\zeta(k)}\log x \right)\\
\end{equation}
For $n=0$ case only, we assume the notation $\gamma^{M,k}=\gamma^{M,k}_0$.

In Table 1, we compute the first $10$ constants for  $k=2$ to $11$ and observe that its value is approaching the Euler's constant as $k\to\infty$, hence one has

\begin{equation}\label{eq:1}
\lim_{k\to\infty}\gamma^{M,k} =\gamma
\end{equation}
the reason being is that the k-free Dirichlet series approaches the harmonic series
\begin{equation}\label{eq:1}
\lim_{k\to\infty}\sum_{n\leq x}^{}\frac{\mu^{(k)}(n)}{n} \to \sum_{n\leq x}^{}\frac{1}{n}
\end{equation}
as $k\to\infty$, where only the higher order $k$ divisors are missing, and one also has $\frac{1}{\zeta(k)}\to 1$.

\begin{table}[hbt!]
\caption{A high-precision computation k-free analogue constant} 
\centering 
\begin{tabular}{| c | c |} 
\hline 
$\textbf{k}$ & \boldmath$\gamma^{M,k}$ \\ [0.5ex] 
\hline 
$2$ & 1.043894515711938297404563438509   \\
\hline
$3$ & 0.891541580568652696893988105635   \\
\hline
$4$ & 0.768619483359480556610282699590   \\
\hline
$5$ & 0.689533694428236052333652469748    \\
\hline
$6$ & 0.641881881833269966325319431727    \\
\hline
$7$ & 0.613974332885360734840225581982    \\
\hline
$8$ & 0.597899163505840636385864829837   \\
\hline
$9$ & 0.588751517287173610371480648659  \\
\hline
$10$ &0.583598636665716799342777885260  \\
\hline
$11$ &0.580722048488931420177890261674   \\
\hline
\end{tabular}
\label{table:nonlin} 
\end{table}

We next consider an alternative representation for these coefficients. If we take the Laurent series expansion of zeta about $s=1$ of

\begin{equation}\label{eq:1}
\zeta(s)=\frac{1}{s-1}+\gamma+O(|s-1|)
\end{equation}
to the $0^{th}$ order, and another expansion of

\begin{equation}\label{eq:1}
\frac{1}{\zeta(ks)}=\frac{1}{\zeta(k)}-k\frac{\zeta^{\prime}(k)}{\zeta(k)^2}(s-1)+O(|s-1|^2)
\end{equation}
to the $1^{st}$ order, then by multiplying them together yields

\begin{equation}\label{eq:1}
\frac{\zeta(s)}{\zeta(ks)}=\frac{1}{\zeta(k)(s-1)}+\frac{\gamma}{\zeta(k)}-k\frac{\zeta^{\prime}(k)}{\zeta(k)^2}+O(|s-1|).
\end{equation}
As a result, the $0^{th}$ order coefficient in (18) is collected in the limit $s\to 1$ as

\begin{equation}\label{eq:1}
\gamma^{M,k} =\frac{\gamma}{\zeta(k)}-k\frac{\zeta^{\prime}(k)}{\zeta(k)^2}
\end{equation}
gives a better closed-form formula in terms of other known constants. And the value for the zeta derivative $\zeta'(k)$ can be computed by

\begin{equation}\label{eq:1}
\zeta^{\prime}(k)=-\sum_{n=1}^{\infty}\frac{\log(n)}{n^k}.
\end{equation}

In Tables 2-5, we compute the first $10$ higher order $\gamma^{M,k}_n$ coefficients to high-precision ($30$ decimal places) by the $n^{th}$ differentiation of the (rhs) of (4) by subtracting the pole term in the limit as

\begin{equation}\label{eq:1}
\gamma^M_{n,k} =\frac{d^n}{ds^n}\Bigg[\frac{\zeta(s)}{\zeta(ks)}-\frac{1}{\zeta(k)(s-1)}\Bigg]_{s\to 1}
\end{equation}
which is easily possible to do in software packages such as Mathematica or Pari/GP [4].

\renewcommand{\figurename}{Table} 
\begin{figure}[hb]
 \setcounter{figure}{1} 
    \begin{minipage}{0.45\textwidth}
        \centering
        \caption{}
       \begin{tabular}{| c | c |} 
\hline 
$\textbf{n}$ & \boldmath$\gamma_n^{M,2}$ \\ [0.5ex] 
\hline 
$0$ & 1.043894515711938297404563438509   \\
\hline
$1$ &-0.236152886477122974860578286060  \\
\hline
$2$ & 0.319384120408014249249465207074  \\
\hline
$3$ &-0.501294458741649566645935631332   \\
\hline
$4$ & 1.010739722784850417039579626049   \\
\hline
$5$ & -2.544030257932552280334481508980   \\
\hline
$6$ & 7.666100995112318690725728704276  \\
\hline
$7$ &-26.88797470534219199661349019865  \\
\hline
$8$ & 107.6566910334506652692812639473  \\
\hline
$9$ &-484.6934692784684121614213582581  \\
\hline
$10$& 2424.080089640181055133479838894  \\
\hline
\end{tabular}
    \end{minipage}
    \hfill
    \begin{minipage}{0.45\textwidth}
        \centering
        \caption{}
      \begin{tabular}{| c | c |} 
\hline 
$\textbf{n}$ & \boldmath$\gamma_n^{M,3}$\\ [0.5ex] 
\hline 
$0$ & 0.891541580568652696893988105635   \\
\hline
$1$ & -0.245232425193549088584910660338   \\
\hline
$2$ & 0.478859132334302048001657659233   \\
\hline
$3$ & -1.001749939953870715953016999262    \\
\hline
$4$ & 2.605658110908472598918158429694     \\
\hline
$5$ & -8.280181525858129142262381587973    \\
\hline
$6$ & 31.06687418066902031238934836818   \\
\hline
$7$ & -134.2881733346003220009976922544   \\
\hline
$8$ & 657.6022908453450809684814267259   \\
\hline
$9$ & -3601.372879333978339048317645190   \\
\hline
$10$& 21824.75413023318146861252990975   \\
\hline
\end{tabular}
    \end{minipage}

  \vspace{1em} 

    \begin{minipage}{0.45\textwidth}
        \centering
        \caption{}
       \begin{tabular}{| c | c |} 
\hline 
$\textbf{n}$ & \boldmath$\gamma_n^{M,4}$ \\ [0.5ex] 
\hline 
$0$ & 0.768619483359480556610282699590    \\
\hline
$1$ & -0.181272887555709074717349472246    \\
\hline
$2$ & 0.482385593391270339523799620921    \\
\hline
$3$ & -1.217994351127279734636351368008     \\
\hline
$4$ & 3.711024193639538830737880019661     \\
\hline
$5$ & -13.61003278671187974663949651367     \\
\hline
$6$ & 58.48358113969641160388773799724   \\
\hline
$7$ & -287.8121812713812744669879624051    \\
\hline
$8$ & 1596.153337468935239880936537014    \\
\hline
$9$ & -9857.607981656058667385016103321    \\
\hline
$10$& 67147.71869990596717226860751284    \\
\hline
\end{tabular}
    \end{minipage}
    \hfill
    \begin{minipage}{0.45\textwidth}
        \centering
        \caption{}
      \begin{tabular}{| c | c |} 
\hline 
$\textbf{n}$ & \boldmath$\gamma_n^{M,5}$\\ [0.5ex] 
\hline 
$0$ & 0.689533694428236052333652469748    \\
\hline
$1$ & -0.110662327010032813709101796076    \\
\hline
$2$ & 0.415273869629961542853619205015    \\
\hline
$3$ & -1.216297879807816360207457417664     \\
\hline
$4$ & 4.165641604898467635390667588548    \\
\hline
$5$ & -16.90556080242211587404036520880     \\
\hline
$6$ & 79.86174202644288489695500264961    \\
\hline
$7$ & -430.5047021058354385867335267643   \\
\hline
$8$ & 2605.890311019817862275989909315    \\
\hline
$9$ & -17505.48225950875316901409735157    \\
\hline
$10$& 129339.3988398803463345898998687    \\
\hline
\end{tabular}
    \end{minipage}
\end{figure}

\newpage

\section{Proof of Theorem 1}

To show that, let us recap some basic definitions. The k-free counting function is defined by

\begin{equation}\label{eq:1}
Q_k(x)=\sum_{n\leq x}\mu^{(k)}(n)
\end{equation}
for positive integer argument $x>1$, and by the average value

\begin{equation}\label{eq:1}
Q_k(x)=\frac{1}{2}\left[\sum_{n< x}\mu^{(k)}(n)+\sum_{n\leq x}\mu^{(k)}(n)\right]
\end{equation}
for positive real argument $x>0$, so whenever $x$ is integer an average value of left and right side of the step is taken. One has the asymptotic
\begin{equation}\label{eq:1}
Q_k(x)\sim \frac{x}{\zeta(k)}
\end{equation}
and the limit
\begin{equation}\label{eq:1}
\lim_{x\to\infty}\frac{Q_k(x)}{x}=\frac{1}{\zeta(k)}
\end{equation}
which implies that the probability of choosing a k-free integer at random from a set of integers is $\frac{1}{\zeta(k)}$. The k-free counting function can be expressed as

\begin{equation}\label{eq:1}
Q_k(x)=\frac{x}{\zeta(k)}+f_k(x)
\end{equation}
where the bound

\begin{equation}\label{eq:1}
f_k(x)=O(x^{\frac{1}{k}})
\end{equation}
which it suffices to prove by the prime number theorem.  But if assuming (RH), there is a much better bound

\begin{equation}\label{eq:1}
f_k(x)=O(x^{\frac{1}{2k}})
\end{equation}
as shown in [3]. Also, the $m^{th}$ derivatives of (1) are given by
\begin{equation}\label{eq:1}
(-1)^m\left[\frac{\zeta(s)}{\zeta(ks)}\right]^{(m)}=\sum_{n=1}^{\infty}\frac{\mu^{(k)}(n)\log^m(n)}{n^s}
\end{equation}
for $\Re(s)>1$.

So we proceed proving Theorem $1$ by carrying out appropriate Stieljtes integration. For $m\geq 0$ case, we have:
\begin{align}
\sum_{n\leq x}\frac{\mu^{(k)}(n)\log^m(n)}{n} &= \int_{1^{-}}^{x} \frac{\log^m(t)}{t}dQ_k(t)
= \int_{1^{-}}^{x} \frac{\log^m(t)}{t}d\left(\frac{t}{\zeta(k)}\right) + \int_{1^{-}}^{x} \frac{\log^m(t)}{t}df_k(t) \notag \\[1.2em]
&= \frac{1}{\zeta(k)}\int_{1^{-}}^{x} \frac{\log^{m}(t)}{t}\, dt
+ \left[ \frac{\log^m(t)}{t}f_k(t) \right]_{1^{-}}^{x} - \int_{1^{-}}^{x} \frac{d}{dt}\left[\frac{\log^{m}(t)}{t}\right]f_k(t)\, dt \notag \\[1.2em]
&= A(x) + \frac{\log^m(x)}{x}f_k(x)+B_m+\int_{1^{-}}^{x}\frac{\log^m(t)-m\log^{m-1}(t)}{t^2}f_k(t)\, dt
\end{align}

\noindent Here we note there is a subtle detail concerning the lower integration variable where we must consider the value

\begin{equation}\label{eq:1}
Q_k(1^{-}) = 0
\end{equation}
before the step, since $Q_k(1^{+}) = 1$ after the step. As a result, the value

\begin{equation}\label{eq:1}
f_k(1^{-})=-\frac{1}{\zeta(k)}
\end{equation}
is the right solution. We define a constant $B_m$ by

\begin{equation}\label{eq:1}
B_m=
\begin{cases}
\frac{1}{\zeta(k)} & \text{if } m = 0 \\
0   & \text{if } m > 0
\end{cases}
\end{equation}
which only arises for $m=0$ case. The reason is that it is assumed the convention that $0^0=1$, which implies $[\log(1)]^0=1$, thus yielding a nonzero term, otherwise for $m>1$, we have $[\log(1)]^m=0^m=0$.

And this integral identity can be generated as

\begin{equation}\label{eq:1}
\begin{aligned}
A(x)&=\frac{1}{\zeta(k)}\int_{1}^{x} \frac{\log^{m}(t)}{t}\, dt\\[1.2em]
&=\frac{\log^{m+1}(x)}{\zeta(k)(m+1)}
\end{aligned}
\end{equation}
for $m\geq 1$. Then, it is readily seen that the limit

\begin{equation}\label{eq:1}
\sum_{n\leq x}\frac{\mu^{(k)}(n)\log^m(n)}{n}-\frac{\log^{m+1}(x)}{\zeta(k)(m+1)}=(-1)^m \gamma_{k,m}^{M}+O\left(x^{\tfrac{1}{2k}-1}\log^m(x)\right)
\end{equation}
exists as $x\to\infty$ and converges to a constant

\begin{equation}\label{eq:1}
(-1)^{m}\gamma_{k,m}^{M}=B_m+\int_{1^{-}}^{\infty}\frac{\log^m(t)-m\log^{m-1}(t)}{t^2}f_k(t)\, dt
\end{equation}
for $m\geq 0$ by the estimate (27) or (28) for $k\geq 2$. Now the next step is to connect this constant to the series expansion in (4) as follows:

\begin{equation}\label{eq:1}
\begin{aligned}
\frac{\zeta(s)}{\zeta(ks)}
&= \int_{1^{-}}^{\infty}\frac{1}{t^{s}}\, dQ_k(t) = \int_{1^{-}}^{\infty} \frac{1}{t^{s}} d\left(\frac{t}{\zeta(k)}\right)\,  + \int_{1^{-}}^{\infty} \frac{1}{t^{s}} df_k(t) \\[1.2em]
&=\frac{1}{\zeta(k)}\int_{1^{-}}^{\infty} \frac{1}{t^{s}} dt\,  + \Biggl[\frac{f_k(t)}{t^s}\Biggr]_{1^{-}}^\infty
 + s\int_{1^{-}}^\infty t^{-s-1}f_k(t)\, dt \\[1.2em]
 &=\frac{1}{\zeta(k)(s-1)}+\frac{1}{\zeta(k)}+s\int_{1^{-}}^\infty t^{-s-1}f_k(t)\, dt \\[1.2em]
\end{aligned}
\end{equation}
The last integral in (37) converges for $\Re(s)>\tfrac{1}{2k}$, therefore it defines an analytical continuation of (lhs) of (37) to that region.
It now remains to show that this integral admits the series expansion about $s=1$, so we consider

\begin{equation}\label{eq:1}
\int_{1^{-}}^{\infty} t^{-s-1} f_k(t)\, dt
\end{equation}
which converges for $\Re(s)>\tfrac{1}{2k}$ assuming (RH). One then generates the exp-log expansion as
\begin{equation}\label{eq:1}
\int_{1^{-}}^{\infty} t^{-2} e^{-(s-1) \log t}f_k(t) \, dt = \sum_{j=0}^{\infty} c_j \frac{(s-1)^j}{j!}
\end{equation}
where its expansion coefficients can be read off as
\begin{equation}\label{eq:1}
c_j = (-1)^j\int_{1^{-}}^{\infty} \frac{\log^j(t)}{t^{2}}f_k(t)  \, dt
\end{equation}
but such a series is only limited to $|s-1|<1$ domain. However, on writing integral term in (37) in $s-1$ domain, then one needs to consider this form instead

\begin{equation}\label{eq:1}
\begin{aligned}
s\int_{1^{-}}^{\infty} t^{-s-1} f_k(t)\, dt &= (s-1)\int_{1^{-}}^{\infty} t^{-s-1} f_k(t)\, dt+\int_{1^{-}}^{\infty} t^{-s-1} f_k(t)\, dt \\[1.2em]
&= \int_{1^{-}}^{\infty} [(s-1)t^{-s-1} f_k(t)+t^{-s-1} f_k(t)]\, dt \\[1.2em]
\end{aligned}
\end{equation}
where we can generate a new series expansion as
\begin{equation}\label{eq:1}
\begin{aligned}
s\int_{1^{-}}^{\infty} t^{-s-1} f_k(t)\, dt&=\sum_{j=0}^{\infty}(j+1)c_j\frac{(s-1)^{j+1}}{(j+1)!}+\sum_{j=0}^{\infty}c_j\frac{(s-1)^{j}}{j!}\\[1.2em]
&=\sum_{m=1}^{\infty}mc_{m-1}\frac{(s-1)^{m}}{m!}+c_0+\sum_{m=1}^{\infty}c_m\frac{(s-1)^{m}}{m!}\\[1.2em]
&=c_0+\sum_{m=1}^{\infty}\left[mc_{m-1}+c_m\right]\frac{(s-1)^{m}}{m!}\\[1.2em]
&=\sum_{j=0}^{\infty} y_j \frac{(s-1)^j}{j!} \\[1.2em]
\end{aligned}
\end{equation}
by combining the expansions in (39). Here we've shifted index variable $m=j+1$ in first series in (42) and re-labeled $m=j$ in second series in (42). The new expansion coefficients are then

\begin{equation}\label{eq:1}
y_m=\begin{cases}
c_0 & \text{for } m = 0 \\
mc_{m-1}+c_m & \text{for } m \geq 1
\end{cases}
\end{equation}
where they are expressed by $c_n$ coefficients, and by (40) these coefficients can be expressed by the integral

\begin{equation}\label{eq:1}
y_m=(-1)^m\int_{1^{-}}^{\infty}\frac{\log^m(t)-m\log^{m-1}(t)}{t^2}f_k(t)\, dt
\end{equation}
which matches the series in (36). And we note that for $m=0$ case we simplify (44) to get

\begin{equation}\label{eq:1}
y_0=\int_{1^{-}}^{\infty}\frac{f_k(t)}{t^2}\, dt
\end{equation}
whereby adding the constant $B_0$ from (33) (and as it appears in (37)), we get the formula

\begin{equation}\label{eq:1}
\gamma^{M,k}=\frac{1}{\zeta(k)}+\int_{1^{-}}^{\infty}\frac{f_k(t)}{t^2}\, dt
\end{equation}
for $k\geq 2$ case. And for the $m\geq 1$ case we have $\gamma^M_{m,k}=y_m$. Hence this completes the proof.

\section{On analytical continuation of the Dirichlet series}
We saw earlier that the k-free Dirichlet series

\begin{equation}\label{eq:1}
\frac{\zeta(s)}{\zeta(ks)}=\sum_{n=1}^{\infty}\frac{\mu^{(k)}(n)}{n^s}
\end{equation}
is convergent for $\Re(s)>1$. We show this by similar Stieltjes integration, but of the remainder term as follows:

 \begin{equation}\label{eq:1}
\begin{aligned}
\frac{\zeta(s)}{\zeta(ks)}
&= \sum_{n\leq x}^{}\frac{\mu^{(k)}(n)}{n^s}+\int_{x}^{\infty}\frac{1}{t^{s}}\, dQ_k(t) \\[1.2em]
& = \sum_{n\leq x}^{}\frac{\mu^{(k)}(n)}{n^s}+\int_{x}^{\infty} \frac{1}{t^{s}} d\left(\frac{t}{\zeta(k)}\right)\,  + \int_{x}^{\infty} \frac{1}{t^{s}} df_k(t) \\[1.2em]
&= \sum_{n\leq x}^{}\frac{\mu^{(k)}(n)}{n^s}+\frac{1}{\zeta(k)}\int_{x}^{\infty} \frac{1}{t^{s}} dt\,  + \Biggl[\frac{f_k(t)}{t^s}\Biggr]_{x}^\infty
 + s\int_{x}^\infty t^{-s-1}f_k(t)\, dt \\[1.2em]
&=\sum_{n\leq x}^{}\frac{\mu^{(k)}(n)}{n^s}-\frac{x^{1-s}}{\zeta(k)(1-s)}-\frac{f_k(x)}{x^s}+s\int_{x}^\infty t^{-s-1}f_k(t)\, dt \\[1.2em]
&=\sum_{n\leq x}^{}\frac{\mu^{(k)}(n)}{n^s}-\frac{s x^{1-s}}{\zeta(k)(1-s)}-\frac{Q_k(x)}{x^s}+O\Bigg(\dfrac{s}{(s-\tfrac{1}{2k})x^{s-\tfrac{1}{2k}}}\Bigg)\, \\[1.2em]
\end{aligned}
\end{equation}

When extracting in the limit we get

\begin{equation}\label{eq:1}
\frac{\zeta(s)}{\zeta(ks)}=\lim_{x\to\infty}\Bigg\{\sum_{n\leq x}^{}\frac{\mu^{(k)}(n)}{n^s}-\frac{x^{1-s}}{\zeta(k)(1-s)}\Bigg\}
\end{equation}
valid for $\Re(s)>\frac{1}{2k}$ for $k\geq 2$ (except at $s=1$) assuming (RH), where we can drop the higher order terms. But if we use the higher order terms in (48) then the precision of will improve.

Now, when $k\to\infty$, then this formula approaches the well-known formula

\begin{equation}\label{eq:1}
\zeta(s)=\lim_{x\to\infty}\Bigg\{\sum_{n\leq x}^{}\frac{1}{n^s}-\frac{x^{1-s}}{1-s}\Bigg\}
\end{equation}
for the zeta function obtained from Euler-Maclaurin formula, which is convergent for $\Re(s)>0$ (except at the pole $s=1$). So basically the k-free Dirichlet series is sort of like an intermediate step for (50).

%

\section{A special case around the first zero}
In the previous section, we have the analytical continuation formula (49) of the k-free Dirichlet series

\begin{equation}\label{eq:1}
\frac{\zeta(s)}{\zeta(ks)}=\lim_{x\to\infty}\Bigg\{\sum_{n\leq x}^{}\frac{\mu^{(k)}(n)}{n^s}-\frac{x^{1-s}}{\zeta(k)(1-s)}\Bigg\}
\end{equation}
converges for $\Re(s)>\frac{1}{2k}$ (except at $s=1$). This implies that all non-trivial zeros of the zeta function on the critical line are also zeros of (51), and even the zeros close the critical line if such exists. Now, there is one more zero that occurs for this representation at $s=\frac{1}{k}$ due to the pole of $\zeta$ in the denominator, and so, for $k\geq 2$ we obtain a special case asymptotic relation at this zero as

\begin{equation}\label{eq:1}
\sum_{n\leq x}^{} \frac{\mu^{(k)}(n)}{n^{\frac{1}{k}}} \sim \frac{k}{(k-1)\zeta(k)} x^{(1-\frac{1}{k})}
\end{equation}
as $x\to \infty$, or expressing in another way, we can extract the reciprocal zeta (for integer $k\geq 2$) as

\begin{equation}\label{eq:1}
\frac{1}{\zeta(k)}=\lim_{x\to\infty}\left(1-\frac{1}{k}\right)x^{(1-\frac{1}{k})}\sum_{n\leq x}^{}\frac{\mu^{(k)}(n)}{n^{\frac{1}{k}}}.
\end{equation}

We also consider another Taylor series expansion about $s=\frac{1}{k}$ of

\begin{equation}\label{eq:1}
\frac{\zeta(s)}{\zeta(ks)}=k\zeta(\frac{1}{k})(s-\frac{1}{k})+O(|s-\frac{1}{k}|^2)
\end{equation}
and seek to extract the first order zeta value by differentiation of the (rhs) of (51). We obtain the formula

\begin{equation}\label{eq:1}
\zeta(\frac{1}{k})=\lim_{x\to\infty}\Bigg\{-\frac{1}{k}
\sum_{n\leq x}\frac{\mu^{(k)}(n)\log(n)}{n^{\frac{1}{k}}}-\frac{1}{\zeta(k)(k-1)}x^{1-\frac{1}{k}}\left(\frac{k}{k-1}-\log(x)\right)\Bigg\}
\end{equation}
And if we substitute (52) to (55) then we generate another type of formula

\begin{equation}\label{eq:1}
\zeta(\frac{1}{k})=\lim_{x\to\infty}\Bigg\{-\frac{1}{k}\sum_{n\leq x}^{}\frac{\mu^{(k)}(n)(\log(n)+\frac{k}{k-1})}{n^{\frac{1}{k}}}+\frac{1}{\zeta(k)(k-1)}x^{1-\frac{1}{k}}\log(x)\Bigg\}
\end{equation}

We will next compute these formulas and verify numerically. First, for the case $k=2$, this formula simplifies to

\begin{equation}\label{eq:1}
\zeta(\frac{1}{2})=\lim_{x\to\infty}\Bigg\{-\frac{1}{2}\sum_{n\leq x}^{}\frac{\mu^{(2)}(n)(\log(n)+2)}{\sqrt{n}}+\frac{6}{\pi^2}\sqrt{x}\log(x)\Bigg\}
\end{equation}
and in Fig.~1, we plot it as a function of limit variable $x$ from $x=1$ to $10^9$ in increments of $1$ on a logarithmic scale on the x-axis to capture all points, thus, every single point is captured in the plot. As a result, we find that this series is a rather noisy and fluctuating series, and we note how it is oscillating about the value $\zeta(\frac{1}{2})=-1.46035450\ldots$, with a gradual decrease in amplitude as $x$ increase (as we shown in [2]).
And similarly, for $k=3$ case we generate the formula

\begin{equation}\label{eq:1}
\zeta(\frac{1}{3})=\lim_{x\to\infty}\Bigg\{-\frac{1}{3}\sum_{n\leq x}^{}\frac{\mu^{(3)}(n)(\log(n)+\frac{3}{2})}{n^{\frac{1}{3}}}+\frac{1}{2\zeta(3)}x^{\frac{2}{3}}\log(x)\Bigg\}
\end{equation}
and in Fig.~2, we plot it as a function of variable x, and show that is oscillating about the value $\zeta(\tfrac{1}{3})= -0.97336024\ldots$ And for $k=4$ case we generate
\begin{equation}\label{eq:1}
\zeta(\frac{1}{4})=\lim_{x\to\infty}\Bigg\{-\frac{1}{4}\sum_{n\leq x}^{}\frac{\mu^{(4)}(n)(\log(n)+\frac{4}{3})}{n^{\frac{1}{k}}}+\frac{1}{3\zeta(4)}x^{\frac{3}{4}}\log(x)\Bigg\}
\end{equation}
where we plot it in Fig.~3, and $\zeta(\tfrac{1}{4})=-0.81327840\ldots$ And for $k=5$ we generate

\begin{equation}\label{eq:1}
\zeta(\frac{1}{5})=\lim_{x\to\infty}\Bigg\{-\frac{1}{5}\sum_{n\leq x}^{}\frac{\mu^{(5)}(n)(\log(n)+\frac{5}{4})}{n^{\frac{1}{5}}}+\frac{1}{4\zeta(5)}x^{\frac{4}{5}}\log(x)\Bigg\}
\end{equation}
where we plot it in Fig.~4, and $\zeta(\tfrac{1}{5})= -0.73392092\ldots$ In all these cases, we observe that these sequences are oscillating about the correct value, and that these fluctuations decay, but slowly.

\renewcommand{\figurename}{Figure}

\begin{figure}[h]
  \renewcommand{\thefigure}{1}%
  \centering
  \includegraphics[width=110mm]{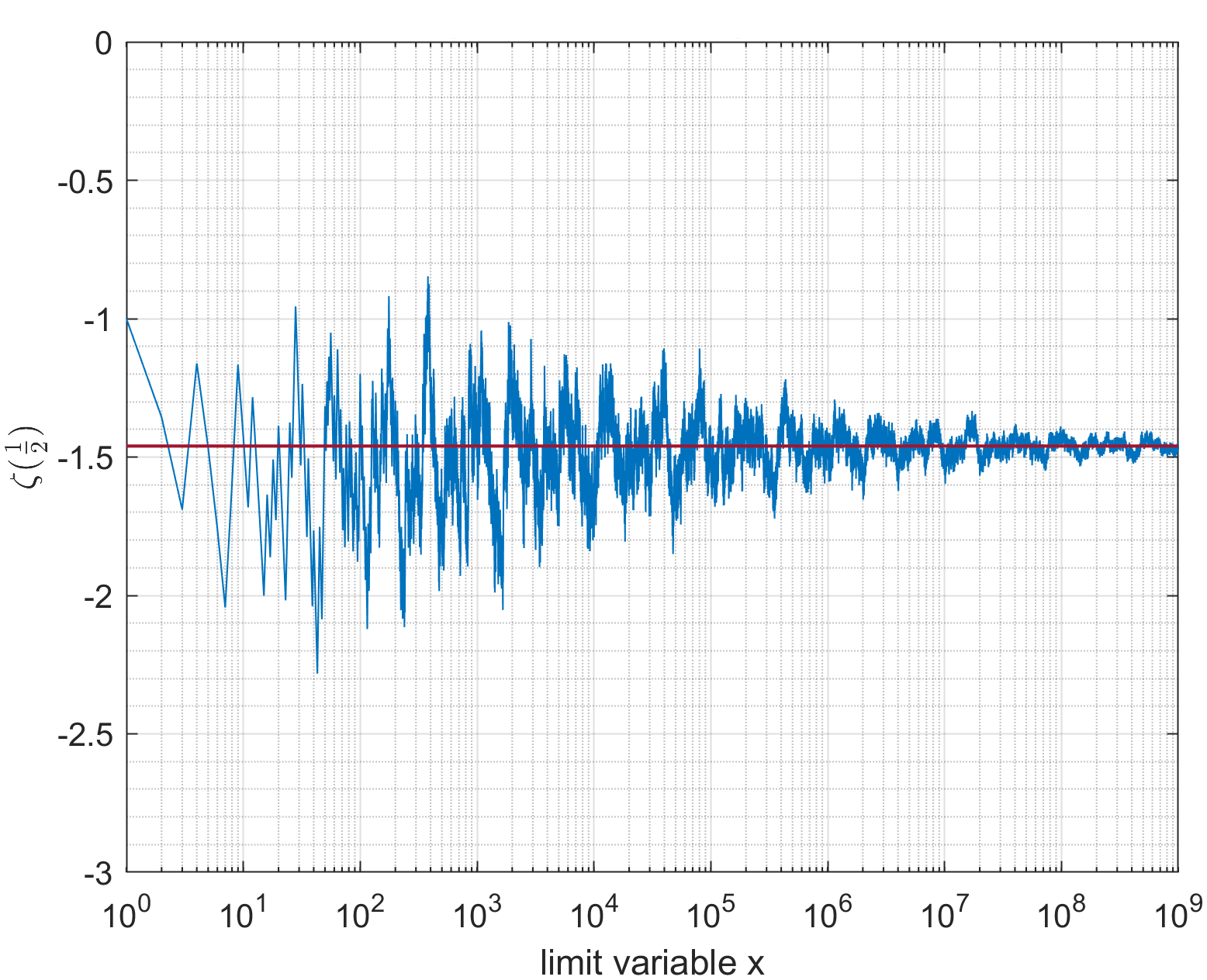}\\
  \caption{A plot of equation (57) for $\zeta(\frac{1}{2})$ value $(k=2)$ as a function of limit variable $x$}\label{1}
  \renewcommand{\thefigure}{2}%
  \centering
  \includegraphics[width=110mm]{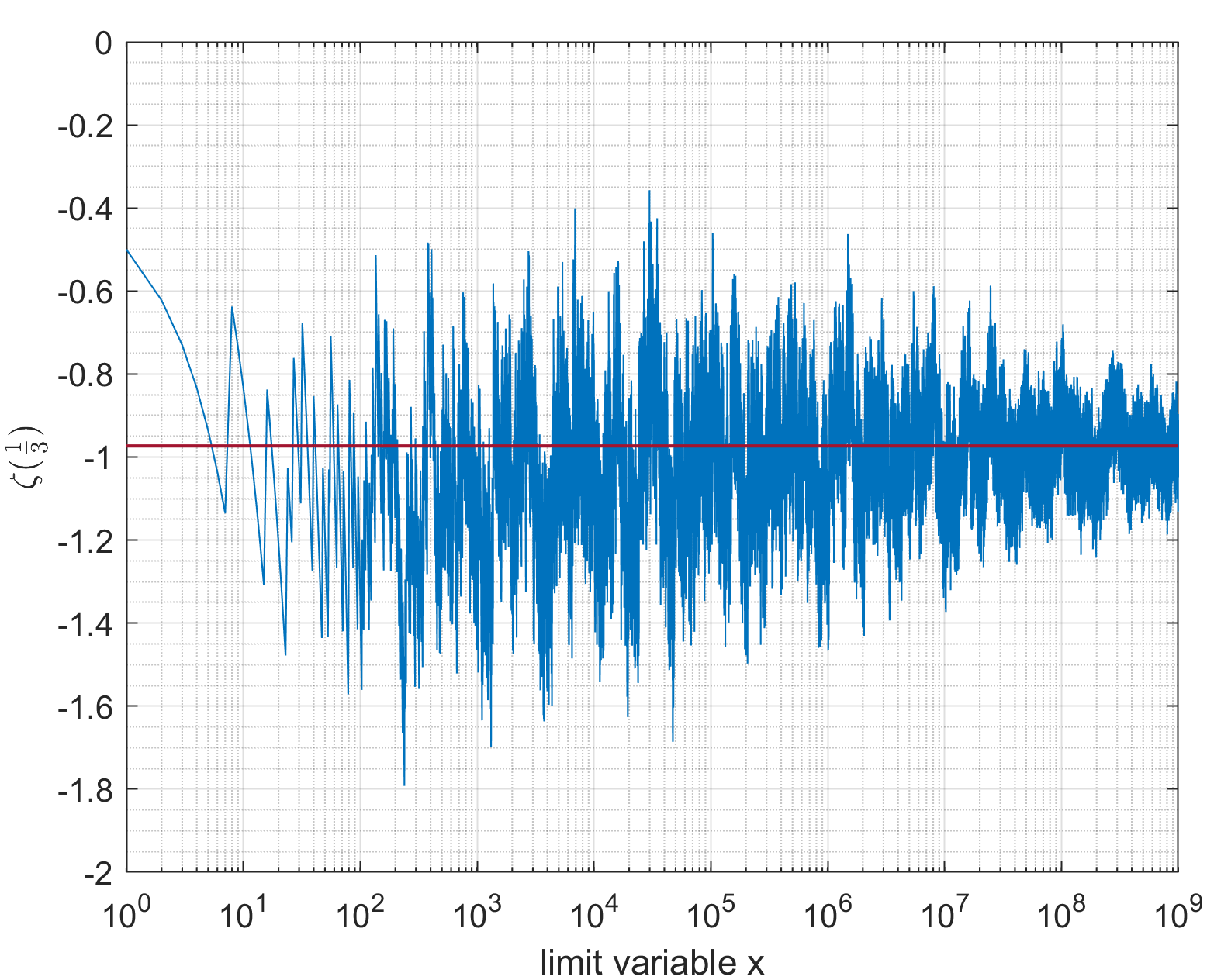}\\
  \caption{A plot of equation (58) for $\zeta(\frac{1}{3})$ value $(k=3)$ as a function of limit variable $x$}\label{1}
\end{figure}

\begin{figure}[htb!]
  \centering
  \renewcommand{\thefigure}{3}%
  \includegraphics[width=110mm]{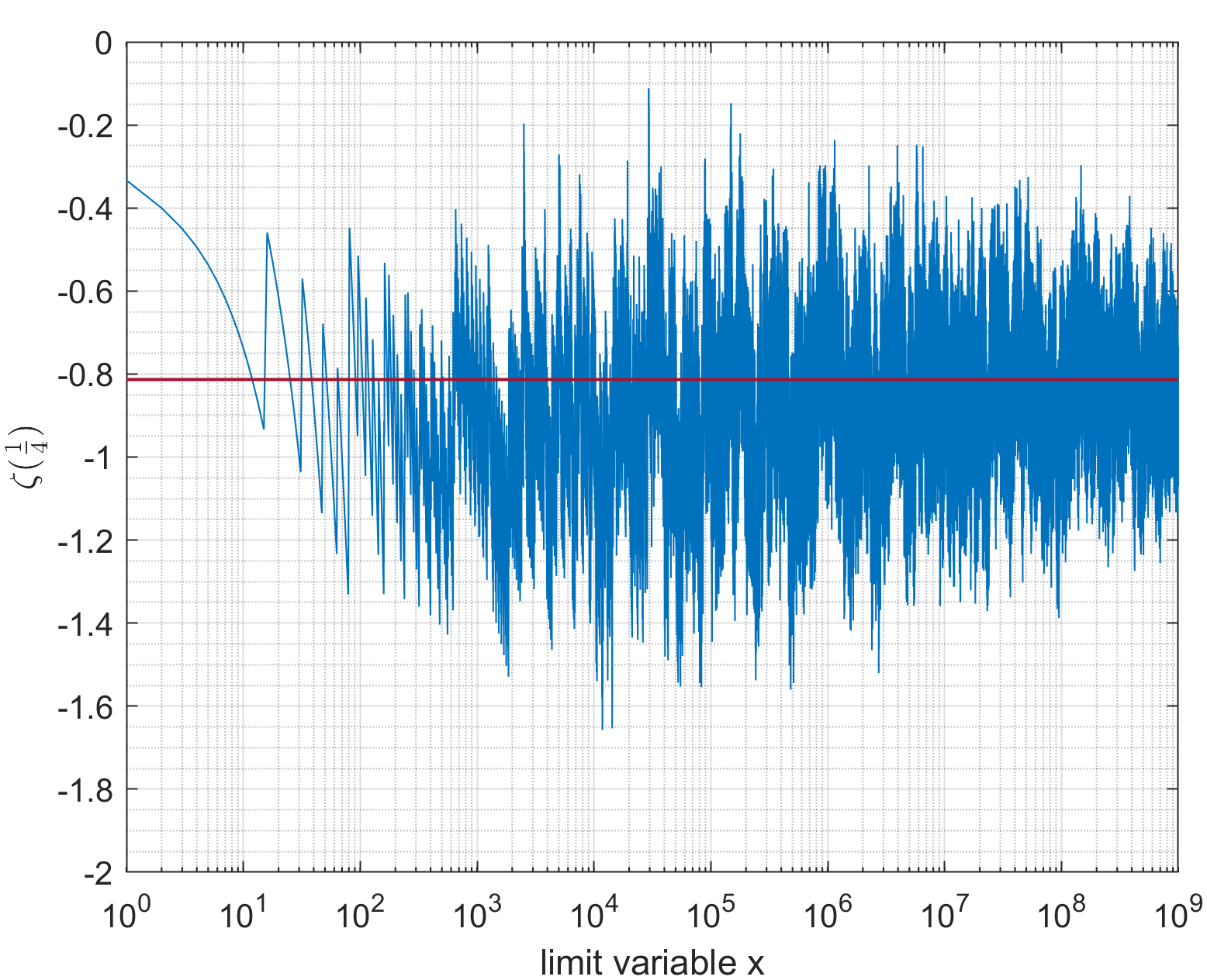}\\
  \caption{A plot of equation (59) for $\zeta(\frac{1}{4})$ value $(k=4)$ as a function of limit variable $x$}\label{1}

   \centering
   \renewcommand{\thefigure}{4}%
  \includegraphics[width=110mm]{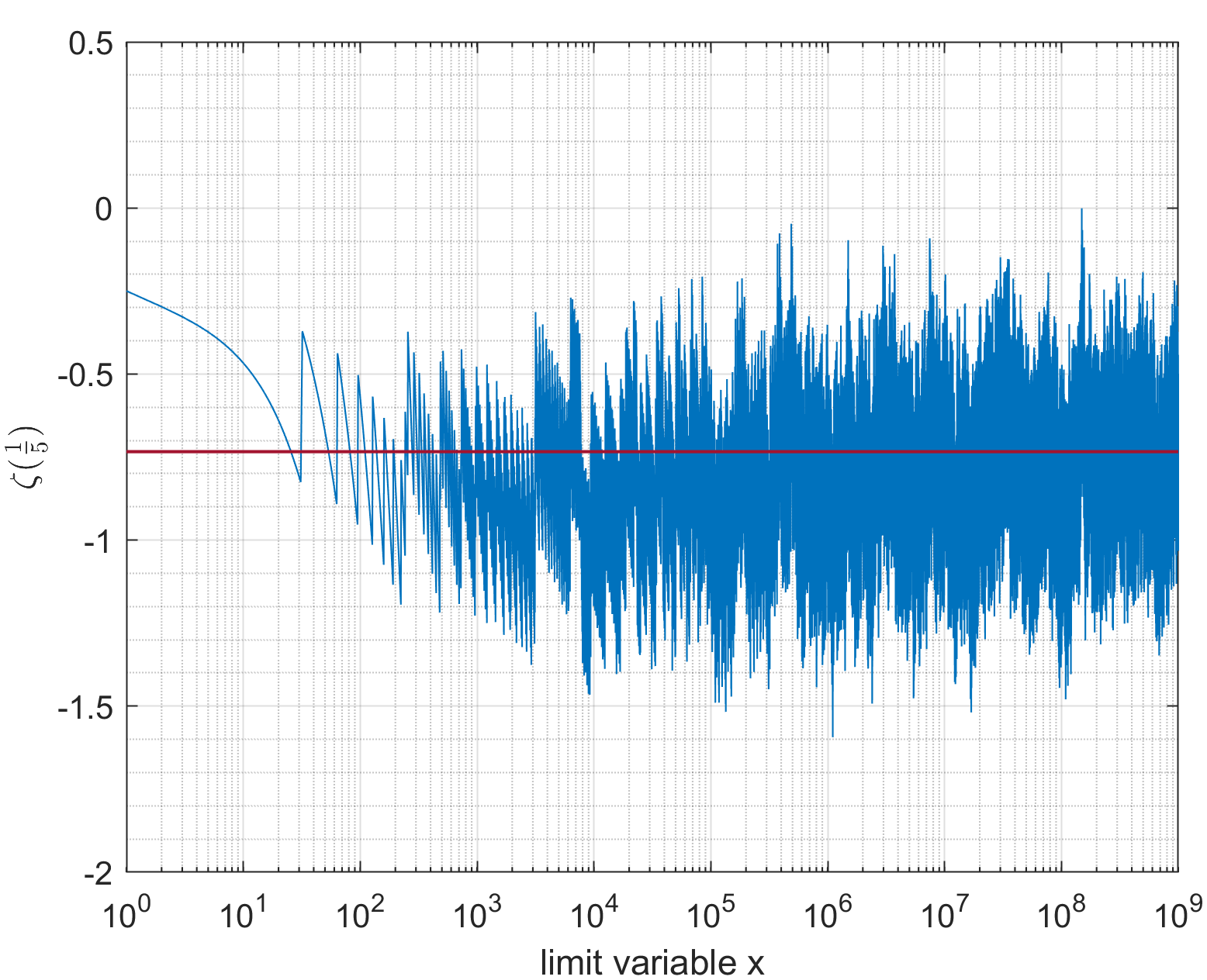}\\
  \caption{A plot of equation (60) for $\zeta(\frac{1}{5})$ value $(k=5)$ as a function of limit variable $x$}\label{1}
\end{figure}


\begin{thebibliography}{9}

\bibitem{latexcompanion}
M. Abramowitz and I. A. Stegun. \textit{Handbook of Mathematical Functions with
Formulas, Graphs, and Mathematical Tables}.
Dover Publications, ninth printing, New York, (1964).

\bibitem{latexcompanion}
A. Kawalec. On the series expansion of a square-free zeta series.math.NT/arXiv:2312.16811, (Dec. 2023).

\bibitem{latexcompanion}
M. Mossinghoff, T. Oliveira e Silva, T. Trudgian. \textit{The distribution of k-free numbers}.math.NT/arXiv:1912.04972, (Jun. 2020).

\bibitem{latexcompanion}
The PARI~Group, PARI/GP version \texttt{2.11.4}, Univ. Bordeaux, (2019).

\bibitem{latexcompanion}
M. Tanaka. \textit{Experiments Concerning the Distribution of Squarefree Numbers}. Proc. Japan Acad. \textbf{55}, Ser. A (1979).

\bibitem{latexcompanion}
Wikipedia \textit{Square-free integer}
\url{https://en.wikipedia.org/wiki/Square-free_integer}, {Accessed: Dec 2025}.

\bibitem{latexcompanion}
M. Wolf, \textit{Numerical determination of a certain mathematical constant related to the M\"obius function}, Computational Methods in Science and Techonology \textbf{29}(1-2),17-20 (2023).

\end{thebibliography}
\end{document}